\documentclass{ifacconf}
\usepackage{amsmath}
\usepackage{graphicx}      
\usepackage{natbib}        
\usepackage[all]{xy}
\usepackage{amssymb,amsfonts,ulem}
\usepackage{stackengine}
\usepackage{scalerel,stackrel}
\usepackage{color}

\newtheorem{theorem}{Theorem}

\newtheorem{definition}[theorem]{Definition}
\newtheorem{example}{Example}
\newtheorem{remark}{Remark}

\begin{document}
\begin{frontmatter}

\title{From retraction maps to symplectic-momentum numerical integrators} 

\thanks[footnoteinfo]{The authors acknowledge financial support from the Spanish Ministry of Science and Innovation under grants PID2022-137909NB-C21, PID2022-137909NB-C22, RED2022-134301-TD, the Severo Ochoa Programme for Centres of Excellence in R\&D (CEX2019-000904-S) and  BBVA Foundation via the project “Mathematical optimization for a more efficient, safer and decarbonized maritime transport”.}

\author[First]{Mar\'ia\ Barbero-Li\~n\'an} 
\author[Second]{Juan Carlos Marrero} 
\author[Third]{David\ Mart\'{\i}n de Diego}

\address[First]{Departamento de Matem\'atica Aplicada, Universidad Polit\'ecnica de Madrid, Av. Juan de Herrera 4, 28040 Madrid, Spain.}
\address[Second]{ULL-CSIC Geometría Diferencial y Mecánica Geométrica, Departamento de Matemáticas Estadística e Investigación Operativa and Instituto de Matemáticas y Aplicaciones (IMAULL), University of La Laguna, San Cristóbal de La Laguna, Spain.}
\address[Third]{Instituto de Ciencias Matem\'aticas (CSIC-UAM-UC3M-UCM), C/Nicol\'as Cabrera 13-15, 28049 Madrid, Spain.}

\begin{abstract}                
Retraction maps have been generalized to discretization maps in~\citep{21MBLDMdD}. Discretization maps are used to systematically derive numerical integrators that preserve the symplectic structure, as well as the discrete momemtum map under the assumption of symmetric preservation for the discretization map. The procedure described here gives a geometrical construction that can be easily adapted to discretize dynamics on more general structures and open the door to reduction processes.  
\end{abstract}

\begin{keyword}
 Symplectic geometry, geometric integrators, symmetry preservation
\end{keyword}

\end{frontmatter}

\section{Introduction}
The study of Hamiltonian systems is widely studied in the context of geometric integration, see for instance \cite{serna,hairer,Feng,blanes}. The main objective is to construct methods preserving the symplecticity.  
Applying backward error analysis, it is possible to show that the modified equations of a symplectic numerical method are also Hamiltonian. This property gives rise to a good energy behaviour, with essentially no accumulation  of error in time for the energy value. 
Apart from the energy, a Hamiltonian system typically has more first integral or constants of motion (as for instance, linear or angular momentum...). Those preservation properties are frequently related with the symmetry invariance of the original Hamiltonian system \citep{AM87}.   

In this paper we focus on symplectic integrators for Hamiltonian systems which also preserve the symmetry of the continuous system. The symmetry condition is introduced in the picture by only imposing the symmetry preservation property to the discretization map. We show that under that condition the symplectic method automatically preserves the associated momentum map. 

\section{Retraction maps and discretization maps}\label{section1}
Retraction maps on manifolds generalizes the linear approximation on Euclidean spaces of a trajectory given a point and a velocity. Let $TM$ be the tangent bundle of a finite-dimensional manifold $M$, a retraction map on a manifold $M$ as defined in~\cite{AbMaSeBookRetraction} is a smooth map 
$R\colon U\subseteq TM \rightarrow M$, where $U$ is an open subset containing the zero section $0_x$  of the tangent bundle at the point $x$ in $M$
such that the restriction map $R_x=R_{|T_xM}\colon T_xM\rightarrow M$ satisfies
\begin{enumerate}
	\item $R_x(0_x)=x$ for all $x\in M$,

	\item and the local rigidity condition
	${\rm D}R_x(0_x)=T_{0_x}R_x={\rm Id}_{T_xM}$ with the identification $T_{0_x}T_xM\simeq T_xM$.
\end{enumerate}
Those retraction maps are used  for optimization theory on Riemannian manifolds replacing the function on $M$ to be optimized by a function on $TM$ that is the Riemmanian gradient of the function to be optimized. The replacement is possible because of the local rigidity condition. 

\begin{example}
If  $(M, g)$ is a Riemannian manifold, then the exponential map $\hbox{exp}^g: U\subset TM\rightarrow M$ is a typical example of retraction map: 
$
\hbox{exp}^g_x(v_x)=\gamma_{v_x}(1), 
$
where $\gamma_{v_x}$ is the unique  Riemannian geodesic (see ~\cite{doCarmo}) satisfying $\gamma_{v_x}(0)=x$ and $\gamma'_{v_x}(0)=v_x$. 
\end{example}

In \cite{21MBLDMdD}  it is defined a generalization of a retraction map called discretization map. It is a map $R_d\colon U \subset TM  \rightarrow M\times M$, where  $U$ is an open neighbourhood of the zero section of $TM$, 
\begin{eqnarray*}
	R_d\colon U \subset TM & \longrightarrow & M\times M\\
	v_x & \longmapsto & (R^1(v_x),R^2(v_x))\, .
\end{eqnarray*}
Discretization maps satisfy the following properties:
\begin{enumerate}
	\item $R_d(0_x)=(x,x)$, for all $x\in M$.
	\item $T_{0_x}R^2_x-T_{0_x}R^1_x \colon T_{0_x}T_xM\simeq T_xM \rightarrow T_xM$ is equal to the identity map on $T_xM$ for any $x$ in $M$.
\end{enumerate}
As an immediate consequence, the discretization map $R_d$ is a local diffeomorphism.

\begin{example}
Examples of discretization maps on Euclidean vector spaces  associated with well-known numerical methods are the following ones:
\begin{itemize}
	\item Explicit Euler method:  $R_d(x,v)=(x,x+v).$
	\item Midpoint rule:  $R_d(x,v)=\left( x-\dfrac{v}{2}, x+\dfrac{v}{2}\right).$
		\item $\theta$-methods with $\theta\in [0,1]$:  $$R_d(x,v)=\left( x-\theta \, v, x+ (1-\theta)\, v\right).$$
\end{itemize}
On a Lie group $G$ the discretization  map $R_d: TG\rightarrow G\times G$ associated with the midpoint rule is:
\[R_d(v_g)=\left(g\,\hbox{exp}\left(-\frac{1}{2}g^{-1}v_g\right), g\,\hbox{exp}\left(\frac{1}{2}g^{-1}v_g\right)\right)\, ,\]
where $v_g\in T_gG$ and $\hbox{exp}:{\mathfrak g}\rightarrow G$ is the exponential map from the Lie algebra ${\mathfrak g}$ to $G$.  
\end{example}

An interesting observation deals with the underlying geometry that can be associated with the above-mentioned notion of discretization maps. From~\cite{Mackenzie} it is known that $Q\times Q$ has a Lie groupoid structure and $TQ$ is the corresponding Lie algebroid. As a future research line, we will prove that the Lie groupoid is the geometric structure that makes possible to extend the notion of discretization maps so that geometric integrators can be defined for a wider family of mechanical systems. Many of the properties satisfied by the above mentioned discretization maps have to deal with the algebraic structure associated with the Lie groupoid (set of identities, inversion map...). This topic will be studied in a forthcoming paper \citep{2024BarMarMar}.

\section{Symplectic-momentum integrators}\label{Sec:Motivation}

In this section we assume that the continuous Hamiltonian  system is invariant by the action of a Lie group. If the discretization map is symmetry preserving (see Definition  \ref{def:symmetry-preserving}), then the discrete flow obtained by the discretization map preserves the momentum map of the original continuous system as proved in Theorem~\ref{Theorem:symmetry}. 
Finally, once a Lie group action appears in the picture, it immediately raises the question of how the original system can be reduced and how a discretization map is defined for the reduced space. Thus, another future research line is to describe how a symmetry preserving discretization map is reduced by the action of a Lie group because that is the seed to obtain geometric integrators preserving the associated reduced structure (for instance, Lie-Poisson brackets) as will be shown in~\cite{2024BarMarMar}.

Consider a free and proper (left)-action $\Phi$ of a Lie group $G$ on the manifold $Q$: $\Phi: G\times  Q\rightarrow Q$. For all $q\in Q$ and $g\in G$, 
$\Phi_{g}: Q\rightarrow Q$ is the map defined by $\Phi_{g}(q)=\Phi (g, q)\equiv g\, q$. 
Such an action satisfies the following properties:
\begin{itemize}
	\item $\Phi(e, q)=q$,
	\item $\Phi(g_1, \Phi(g_2, q))=\Phi(g_1g_2, q)$,
\end{itemize}
for all $g_1, \, g_2\in G$, $q\in Q$ and where $e$ denotes the identity element of $G$. Moreover, $\Phi_{g}: Q\rightarrow Q$ is a diffeomorphism for all $g\in G$, see \cite{AM87}. 

For every element $\xi$ in the Lie algebra ${\mathfrak g}$ of $G$, the following vector field on $Q$ is defined: 
\[
\xi_Q(q)=\frac{d}{dt}\Big|_{t=0} \Phi_{{\rm exp}(t\xi)}(q),
\]
that it is called the {\bf infinitesimal generator} of the action  $\Phi$ corresponding to $\xi$. 

\begin{definition}\label{def:symmetry-preserving}
A discretization map $R_d\colon U\subseteq TQ  \rightarrow Q\times Q$ is {\bf symmetry-preserving} if it satisfies that: 
\[
R_d(T_q\Phi_g(v_q))=(\Phi_g\times \Phi_g) (R_d(v_q)) 
\]
for all $g\in G$ and  $v_q\in U\subseteq TQ$.
\end{definition}
Definition~\ref{def:symmetry-preserving} is equivalent to the commutativity of the following diagram:
\begin{equation*}
\xymatrix{
	 TQ \ar[d]<2pt>_{T\Phi_{g}} \ar[rr]^{R_d} && Q\times Q \ar[d]^{\Phi_g\times \Phi_g}   \\
	 TQ  \ar[rr]^{R_d} && Q\times Q  }
\end{equation*}

As $R_d(v_q)= (R^1_d(v_q),R^2_d(v_q))$, a symmetry-preserving discretization map infinitesimally satisfies
\[
T_{v_q}R_d^i(\xi_Q^C(v_q))=\xi_Q (R_d^i(v_q)), \qquad i=1, 2, 
\]
where $\xi^C_Q$
is the complete lift of the vector field $\xi_Q$ on $Q$ to the tangent bundle. Observe that  $\xi^C_Q$ is the infinitesimal generator $\xi_{TQ}$ for the tangent lift of the action corresponding to $\xi$.
In other words, $\xi^C_Q=\xi_{TQ}$. It is known that if $\{\phi_t\}$ is the local flow of $\xi_Q$ then $\{T\phi_t\}$ is the local flow of $\xi_Q^C$. 
\begin{remark}\label{remark-2}
For any discretization map  $R_d:U\subset TQ\rightarrow Q\times Q$ it is possible to define the adjoint discretization 
$R^*_d: \overline{U}\subset TQ\rightarrow Q\times Q$ with $\overline{U}=\{ v_q\in TQ\; |\, -v_q\in U\}$  defined by 
$
R^*_d (v_q)=\left(I_Q \circ R_d\right)(-v_q)
$, where $I_Q(q, q')=(q', q)$ is the inversion map. Such a map
is also a discretization map.  Thus, if $R_d$ is a symmetry preserving discretization map, it is easy to check that $R_d^*$ is also symmetry preserving: 
\begin{align*}
R^*_d(T_q\Phi_g(v_q))&=\left(I_Q \circ R_d\right)(T_q\Phi_g(-v_q))\\
&= I_Q((\Phi_g\times \Phi_g) (R_d(-v_q))) \\
&= (\Phi_g\times \Phi_g)(I_Q\circ R_d)(-v_q) \\
&= (\Phi_g\times \Phi_g)(R^*_d (v_q))
\end{align*}
\end{remark}

Denote by $T^*\Phi: G\times T^*Q\rightarrow T^*Q$ the corresponding {\bf cotangent lift action} given by
\[
(T^*\Phi)_{g}(\alpha_q)=T^*_{\Phi_{g}(q)}\Phi_{g^{-1}}(\alpha_q),\quad g\in G, \quad \alpha_q\in T^*_qQ\, .
\]	
As it is a cotangent lift, it consequently satisfies that 
\[
\pi_Q\circ (T^*\Phi)_g=\Phi_g\circ \pi_Q\qquad \hbox{and} \qquad (T^*\Phi)^*_g\omega_Q=\omega_Q\,,
\]
where $\omega_Q$ is the canonical symplectic form on $T^*Q$.

	 Let ${\mathfrak g}^*$ be the dual Lie algebra of $G$, there exists a $G$-equivariant {\bf momentum map} $J: T^*Q\rightarrow {\mathfrak g}^*$  defined by
	 \begin{equation}\label{eq:momentum}
	 J(\alpha_q)(\xi)=\langle \alpha_q, \xi_Q(q)\rangle\, 
	\end{equation}
	for any $\xi\in {\mathfrak g}$ and $\alpha_q\in T^*Q$. It is possible to consider the function $J_{\xi}: T^*Q\rightarrow {\mathbb R}$ given by $J_{\xi}(\alpha_q)= J(\alpha_q)(\xi)$. The corresponding Hamiltonian vector field $X_{J_\xi}$, that is,    $i_{X_{J_{\xi}}}\omega_Q={\rm d}{J_{\xi}}$, is the infinitesimal generator for the cotangent lift of the action corresponding to $\xi$:
	$\xi_{T^*Q}=X_{J_{\xi}}$.  In fact, $\xi_{T^*Q}$ is just the complete lift $\xi_Q^{C^*}$ of $\xi_Q$ to $TQ$ (for more details, see \cite{AM87}).

	 
  Now,  consider  the canonical antisymplectomorphism $ {\mathcal I}_{TQ}$ between the symplectic manifolds $(T^*T^*Q, \omega_{T^*Q})$ and $(T^*TQ, \omega_{TQ})$ ~\citep{XuMac} that is locally given by: 
	 \begin{equation}\label{eq:I_TQlocal}
	   	 {\mathcal I}_{TQ}(q,p, \mu_q, \mu_p)=(q, \mu_p, -\mu_q, p) \, .  
	 \end{equation}

	 It is simple to prove that 
	 \begin{equation}\label{eq:pairingI_TQ} 
	 \hspace{-0.35cm}\langle  {\mathcal I}_{TQ}(\Lambda_{\alpha_q}), X^C( \pi_{TQ}( {\mathcal I}_{TQ}(\Lambda_{\alpha_q})))\rangle
	 =-\langle \Lambda_{\alpha_q}, X^{C^*}(\alpha_q)\rangle, 
	 \end{equation}
	for all $\Lambda_{\alpha_q}\in T^*_{\alpha_q}T^*Q$,where $X^C$ and $X^{C^*}$ are the complete lifts of a vector field $X$ on $Q$ to $TQ$ and $T^*Q$, respectively. The above equality is satisfied, in particular, for the infinitesimal generators $\xi_{TQ}=\xi_Q^C$ and $\xi_{T^*Q}=\xi_Q^{C^*}=X_{J_{\xi}}$. 
	
	 On the other hand, any discretization map $R_d: TQ\rightarrow Q\times Q$ can be cotangently lifted to $T^*R_d\colon T^*(TQ) \rightarrow T^*(Q\times Q)$. The map $T^*R_d$ is a symplectomorphim between $(T^*TQ, \omega_{TQ})$ and $(T^*(Q\times Q), \omega_{Q\times Q})$~\citep{AM87}.
 \begin{equation*}
\xymatrix{
	T^*TQ\ar[rr]^{T^*R_d} \ar[d]<2pt>_{\pi_{TQ}} && T^*(Q\times Q)\ar[d]^{\pi_{Q\times Q}}\\
	TQ  \ar[rr]^{R_d} && Q\times Q}
\end{equation*}

From the Hamiltonian point of view, the antisymplectomorphism
	 $$(T^*R_d\circ  {\mathcal I}_{TQ}): (T^*T^*Q, \omega_{T^*Q})\longrightarrow (T^*(Q\times Q), \omega_{Q\times Q})$$ is needed to describe the discrete dynamics:
	 \begin{equation*}
\xymatrix{
	T^*T^*Q\ar[rrr]^{(T^*R_d\circ  {\mathcal I}_{TQ})} \ar[d]<2pt>_{\pi_{TQ}\circ   {\mathcal I}_{TQ}} &&& T^*(Q\times Q)\ar[d]^{\pi_{Q\times Q}}\\
	TQ  \ar[rrr]^{R_d} &&& Q\times Q}
\end{equation*}  
	
  \subsection{Symplectic integrators based on discretization maps}
  Given a  Hamiltonian function $H: T^*Q\rightarrow {\mathbb R}$ and its corresponding Hamiltonian vector field
  \[
  i_{X_H}\omega_Q={\rm d} H
  \]
  we 
  define the Lagrangian submanifold ${\rm d} H (T^*Q)=\hbox{Im}\, {\rm d}\,H$  of the symplectic manifold $(T^*T^*Q, \omega_{T^*Q})$ 
 that induces the following Lagrangian submanifold  
	 \begin{equation}\label{eq:DH_method}
	 D^h_H=(T^*R_d\circ  {\mathcal I}_{TQ}) (h \, \hbox{Im}\, {\rm d}\,H)
	 \end{equation}  
  of $ (T^*(Q\times Q), \omega_{Q\times Q})$,  where $h>0$ is the time step size. The above diagram leads to the same set $D^h_H$ as the cotangent lift $\widehat{R_d}: TT^*Q\rightarrow T^*Q\times T^*Q$ of a discretization map introduced in \citep{21MBLDMdD}, but it is obtained by discretizing two different spaces: $T^*T^*Q$ and $TT^*Q$, respectively.

  Note that the manifold
 $T^*(Q\times Q)$ has the structure of symplectic groupoid over $T^*Q$. As mentioned in~\citep{IMMP}, that structure induces a discrete evolution determined by the corresponding  source and target maps, $\alpha$ and $\beta$, respectively: 
	 \begin{eqnarray*}
	 	\alpha (\mu_{q_k}, \mu_{q_{k+1}})&=&-\mu_{q_k}\, ,\\
	 	\beta (\mu_{q_k}, \mu_{q_{k+1}})&=&\mu_{q_{k+1}}\, .
	 \end{eqnarray*}
	 	The  induced discrete dynamics is given by 
	 	\begin{equation*}
	 		(-\mu_{q_k}, \mu_{q_{k+1}})\in D_{H}^h\, ,
	 	\end{equation*}
	and it implicitly determines a symplectic integrator $\Psi_h (\mu_{k})=\mu_{k+1}$ for the Hamiltonian system determined by $H$
 , where $\Psi_h\colon T^*Q \rightarrow T^*Q$ denotes the discrete flow. 
	 \begin{example}
	 From the midpoint discretization map
	 \newline $R_d(x,v)=\left(x-\dfrac{v}{2},x+\dfrac{v}{2}\right)$, we define 
$$T^*R_d (x, v, p_x, p_v)=\left(x-\dfrac{v}{2}, \dfrac{p_x}{2}-p_v; x+\dfrac{v}{2}, \dfrac{p_x}{2}+p_v\right).$$

Therefore,
\[
D^h_H=\left\{(q_k, -p_k; q_{k+1}, p_{k+1})\, \left|\; 
\begin{matrix}
q_k=q-\frac{h}{2}\frac{\partial H}{\partial p}(q,p)\\
p_k=p+\frac{h}{2}\frac{\partial H}{\partial q}(q,p)\\
q_{k+1}=q+\frac{h}{2}\frac{\partial H}{\partial p}(q,p)\\
p_{k+1}=p-\frac{h}{2}\frac{\partial H}{\partial q}(q,p)
\end{matrix}\right.
\right\}\, .
\]
	
	It can be obtained that $q=\frac{q_k+q_{k+1}}{2}$ and $p= \frac{p_k+p_{k+1}}{2}$. Thus, the symplectic method is the following one: 
	\begin{eqnarray*}
	\frac{q_{k+1}-q_k}{h}&=&\frac{\partial H}{\partial p}\left( \frac{q_k+q_{k+1}}{2},\frac{p_k+p_{k+1}}{2}\right)\, ,\\
	\frac{p_{k+1}-p_k}{h}&=&-\frac{\partial H}{\partial q}\left( \frac{q_k+q_{k+1}}{2},\frac{p_k+p_{k+1}}{2}\right)\, .\\
	\end{eqnarray*}

	\end{example}
\begin{remark}
   From the Lagrangian point of view, the construction is even more direct since it is not necessary to use the map ${\mathcal I}_{TQ}$. In fact, a  Lagrangian  function $L: TQ\rightarrow {\mathbb R}$ defines 
the Lagrangian submanifold ${\rm d}\, L (TQ)=\hbox{Im}\, {\rm d}\,L$ of the symplectic manifold $ (T^*TQ, \omega_{TQ})$. Then, the discrete dynamics is given by the Lagrangian submanifold  
\[
	 D^h_L=(T^*R_d) (h \, {\rm d} L(TQ))
	 \]
of $ (T^*(Q\times Q), \omega_{Q\times Q})$, where $h>0$ is the fixed time step size. 
\end{remark}
\subsection{Symmetry preservation}
Under this geometrical framework, it can be obtained the preservation of the momentun map for the discrete flow when a  symmetry preserving retraction map is used to define a geometric integrator (see also \cite{GeMarsden}).
\begin{theorem}\label{Theorem}
    \label{Theorem:symmetry} Let $H\colon T^*Q\rightarrow \mathbb{R}$ be a Hamiltonian function and $\xi\in {\mathfrak g}$. If $J_{\xi}$ is a constant of motion of the Hamiltonian vector field $X_{H}$,  then the method in Equation~\eqref{eq:DH_method} derived from a symmetry-preserving discretization map $R_d$ satisfies that 
\[
J_{\xi}(\mu_{q_k})=J_{\xi}(\mu_{q_{k+1}})
\]	
for all $(-\mu_{q_k}, \mu_{q_{k+1}})\in D^h_{H}$. 
\end{theorem}
\textit{Proof.}
Since $J_{\xi}$ is a constant of the motion of the Hamiltonian vector field $X_H$, it is satisfied $X_{H}(J_{\xi})=0$. Therefore, it holds that $\xi_{T^*Q}(H)=0$.
Observe now that
\begin{eqnarray*}
\langle
T^*R_d\circ {\mathcal I}_{TQ}\circ {\rm d}\, H ,   (\xi_Q , \xi_Q )\rangle
&=&\langle   {\mathcal I}_{TQ}\circ {\rm d}\, H, \xi_{TQ}
\rangle\\
&=&-\langle {\rm d}\, H, \xi_{T^*Q}\rangle=0 
\end{eqnarray*}
because the discretization map $R_d$ is symmetry-preserving and Equation~\eqref{eq:pairingI_TQ} holds.
Using~\eqref{eq:momentum} we deduce that 
\[
J_{\xi}(\mu_{q_k})=J_{\xi}(\mu_{q_{k+1}})
\]	
for all $(-\mu_{q_k}, \mu_{q_{k+1}})\in D_{H}$. 
\hfill $\square$

Note that we are not asking for the Hamiltonian to be symmetry invariant for all the Lie group actions. The hypothesis in Theorem~\ref{Theorem:symmetry} only requires one element $\xi$ in $ {\mathfrak g}$ such that the Hamiltonian function has $J_{\xi}: T^*Q\rightarrow {\mathbb R}$ as a particular constant of motion. Typically, it is simpler to ask for a symmetry preserving discretization map for a bigger set of elements of the Lie group as shows the next example. 

	\begin{example}
	Consider the Lie group $\hbox{GL}(n)$ that acts on ${\mathbb R}^n$ by matrix multiplication on the left, i.e.,
	$\Phi_M (x)= Mx$. 
For $\alpha\in [0, 1]$, the discretization maps   
	\[
	\begin{array}{rrcl}
R_d^{\alpha}:& T{\mathbb R}^n&\longrightarrow& {\mathbb R}^n\times {\mathbb R}^n\\
&(x, v)&\longmapsto &(x+(\alpha-1) v, x+\alpha v)
\end{array}
	\]
	are $\hbox{GL}(n)$-symmetry-preserving since 
	\begin{align*}
	R_d^{\alpha}(Mx, Mv)&=(Mx+(\alpha-1) Mv, Mx+\alpha Mv)\\ &=(M(x+(\alpha-1) v), M(x+\alpha v)) \\
 &=M\, R_d^{\alpha}(x,v)\,.
	\end{align*}
Consider the cotangent lift of this action: 
\[
(T^*\Phi_M)(x, p)=(Mx, (M^{-1})^{T}p)\, ,
\]
with momentum map
\[
J_{\mathfrak  m}(x, p)=x^T {\mathfrak  m} p
\]
for ${\mathfrak  m}\in {\mathfrak g}l(n)$. 
If  $J_{\mathfrak  m}$ is a constant of motion of $H: T^*{\mathbb R}^n\rightarrow {\mathbb R}$, then 
Theorem~\ref{Theorem:symmetry} guarantees that $J_{\mathfrak  m}$ is also a constant of motion of the discrete method obtained using $R^{\alpha}_d$ as in Equation~\eqref{eq:DH_method}.
\end{example}
\subsection{Higher-order symplectic-momentum integrators}
Considering numerical methods based on discretization maps we finally obtain a Lagrangian submanifold $D_H^h$ of $(T^*(Q\times Q), \omega_{Q\times Q})$. We can apply the classical theory of composition of Lagrangian submanifolds (see \cite{guilleminsternberg80}) to increase the order while preserving symplecticity and constants of motion~\citep{yoshida, hairer} associated to some symmetry of the continuous system. For instance, taking
the Lagrangian submanifold obtained by composition: 
\[
D_H^{\gamma_s h}\circ \ldots D_H^{\gamma_1 h}
\]
with 
$
\gamma_1+\ldots+\gamma_s=1$ and $ 
\gamma_1^{p+1}+\ldots+\gamma_s^{p+1}=1
$
we derive a method at least of order $p+1$. Additionally we can use the adjoint discretization maps, as in Remark \ref{remark-2}, to increase the order of the methods. See \cite{21MBLDMdD} for more details. 
\section{Conclusions}
In this paper we have introduced symmetry preserving discretization maps and we have shown their utility to derive geometric preserving integrators. In Theorem \ref{Theorem} we have proved the preservation of some constants of motion associated to the invariant properties of the continuous Hamiltonian. These results will be also useful in optimization theory where the function to be minimized is also $G$-invariant. 

Finally, the results of this paper open the possibility to study the geometrical derivation of numerical integrators for reduced systems preserving the corresponding Lie-Poisson bracket. For that, it will be necessary to extend the notion of retraction and discretization maps to Lie groupoids and Lie algebroids \cite{2024BarMarMar}. 

\bibliography{References}

\begin{thebibliography}{15}
\providecommand{\natexlab}[1]{#1}
\providecommand{\url}[1]{\texttt{#1}}
\providecommand{\urlprefix}{URL }
\expandafter\ifx\csname urlstyle\endcsname\relax
  \providecommand{\doi}[1]{doi:\discretionary{}{}{}#1}\else
  \providecommand{\doi}{doi:\discretionary{}{}{}\begingroup
  \urlstyle{rm}\Url}\fi

\bibitem[{Abraham and Marsden(1987)}]{AM87}
Abraham, R. and Marsden, J. (1987).
\newblock \emph{Foundations of Mechanics}.
\newblock Addison Wesley, second edition.

\bibitem[{Absil et~al.(2008)Absil, Mahony, and
  Sepulchre}]{AbMaSeBookRetraction}
Absil, P.A., Mahony, R., and Sepulchre, R. (2008).
\newblock \emph{Optimization algorithms on matrix manifolds}.
\newblock Princeton University Press, Princeton, NJ.
\newblock With a foreword by Paul Van Dooren.

\bibitem[{Barbero Li\~n\'an et~al.(work in progress)Barbero Li\~n\'an, Marrero,
  and Mart\'in~de Diego}]{2024BarMarMar}
Barbero Li\~n\'an, M., Marrero, J.C., and Mart\'in~de Diego, D. (work in
  progress).
\newblock Retraction maps: A seed of geometric integrators. {P}art {II}:
  {S}ymmetry and reduction, work in progress.

\bibitem[{Barbero Li\~n\'an and Mart\'in~de Diego(2022)}]{21MBLDMdD}
Barbero Li\~n\'an, M. and Mart\'in~de Diego, D. (2022).
\newblock Retraction maps: a seed of geometric integrators.
\newblock \emph{Found. Comput. Math.}
\newblock \doi{10.1007/s10208-022-09571-x}.

\bibitem[{Blanes and Casas(2016)}]{blanes}
Blanes, S. and Casas, F. (2016).
\newblock \emph{A concise introduction to geometric numerical integration}.
\newblock Monographs and Research Notes in Mathematics. CRC Press, Boca Raton,
  FL.

\bibitem[{do~Carmo(1992)}]{doCarmo}
do~Carmo, M.P. (1992).
\newblock \emph{Riemannian geometry}.
\newblock Mathematics: Theory \& Applications. Birkh\"{a}user Boston, Inc.,
  Boston, MA.

\bibitem[{Feng and Qin(2010)}]{Feng}
Feng, K. and Qin, M. (2010).
\newblock \emph{Symplectic geometric algorithms for {H}amiltonian systems}.
\newblock Zhejiang Science and Technology Publishing House, Hangzhou; Springer,
  Heidelberg.

\bibitem[{Guillemin and Sternberg(1980)}]{guilleminsternberg80}
Guillemin, V. and Sternberg, S. (1980).
\newblock The momentum map and collective motion.
\newblock \emph{Ann. of Phys.}, 1278, 220--253.

\bibitem[{Hairer et~al.(2010)Hairer, Lubich, and Wanner}]{hairer}
Hairer, E., Lubich, C., and Wanner, G. (2010).
\newblock \emph{Geometric numerical integration}, volume~31 of \emph{Springer
  Series in Computational Mathematics}.
\newblock Springer, Heidelberg.
\newblock Structure-preserving algorithms for ordinary differential equations,
  Reprint of the second (2006) edition.

\bibitem[{Iglesias-Ponte et~al.(2013)Iglesias-Ponte, Marrero, Mart\'{\i}n~de
  Diego, and Padr\'{o}n}]{IMMP}
Iglesias-Ponte, D., Marrero, J.C., Mart\'{\i}n~de Diego, D., and Padr\'{o}n, E.
  (2013).
\newblock Discrete dynamics in implicit form.
\newblock \emph{Discrete Contin. Dyn. Syst.}, 33(3), 1117--1135.

\bibitem[{Mackenzie and Xu(1994)}]{XuMac}
Mackenzie, K.C.H. and Xu, P. (1994).
\newblock Lie bialgebroids and {P}oisson groupoids.
\newblock \emph{Duke Math. J.}, 73(2), 415--452.

\bibitem[{Mackenzie(2005)}]{Mackenzie}
Mackenzie, K.C.H. (2005).
\newblock \emph{General theory of {L}ie groupoids and {L}ie algebroids}, volume
  213 of \emph{London Mathematical Society Lecture Note Series}.
\newblock Cambridge University Press, Cambridge.

\bibitem[{Sanz-Serna and Calvo(1994)}]{serna}
Sanz-Serna, J.M. and Calvo, M.P. (1994).
\newblock \emph{Numerical {H}amiltonian problems}, volume~7 of \emph{Applied
  Mathematics and Mathematical Computation}.
\newblock Chapman \& Hall, London.

\bibitem[{Yoshida(1990)}]{yoshida}
Yoshida, H. (1990).
\newblock Construction of higher order symplectic integrators.
\newblock \emph{Phys. Lett. A}, 150(5-7), 262--268.

\bibitem[{Zhong and Marsden(1988)}]{GeMarsden}
Zhong, G. and Marsden, J.E. (1988).
\newblock Lie-{P}oisson {H}amilton-{J}acobi theory and {L}ie-{P}oisson
  integrators.
\newblock \emph{Phys. Lett. A}, 133(3), 134--139.

\end{thebibliography}

\end{document}